\documentclass[12pt,english]{amsart}

\usepackage{geometry}
\usepackage{microtype} 

\usepackage{graphicx}   
\usepackage{wrapfig}    
\usepackage{caption}    
\usepackage{subcaption} 
\usepackage{tikz}       

\usepackage{amsmath, amssymb, amsthm, mathtools, esint}
\usepackage{amsfonts} 
\usepackage{enumitem} 

\usepackage[hidelinks, bookmarksdepth=3]{hyperref}
\hypersetup{
    colorlinks=true,
    linkcolor=blue,
    citecolor=blue,
    urlcolor=blue,
}

\usepackage[alphabetic,msc-links,backrefs]{amsrefs}
\usepackage{doi}
\renewcommand{\PrintDOI}[1]{\doi{#1}}

\numberwithin{equation}{section}

\theoremstyle{definition}

\theoremstyle{remark}

\allowdisplaybreaks[4]       

\usepackage{xcolor}

\begin{document}

\title[Mathematics and Sustainability]{Teaching Mathematics and Sustainability \\ (A narration of a rough journey)}

\author{Henrik Shahgholian}
\address{Department of Mathematics, Royal Institute of Technology, Stockholm, Sweden}
\email{henriksh@kth.se}

  \thanks{{\bf Disclaimer: } These notes reflect personal experiences in course development, highlighting challenges and offering insights to inspire similar efforts. }

\makeatletter
\@namedef{subjclassname@2020}{\textup{2020} Mathematics Subject Classification}
\makeatother




\maketitle

\begin{abstract}
This note presents reflections drawn from my recent experiences in teaching a course on mathematics and sustainability, with a particular emphasis on raising awareness of the topic and its broader implications. The lectures were structured to bridge conceptual understanding with practical engagement — combining mathematical modeling, critical discussion, and case studies on real-world sustainability challenges. Much of the actual work was done through interactive problem-solving sessions, group discussions, and short projects where students explored how quantitative reasoning can inform sustainable decision-making. The aim was not only to convey technical tools, but also to cultivate a deeper appreciation of how mathematics can serve as a framework for interpreting and addressing complex societal issues.
\end{abstract}
\section{Introduction}

\subsection{Where it all began}

\begin{wrapfigure}[16]{r}{0.4\textwidth} 
    \centering
    \includegraphics[width=0.4\textwidth]{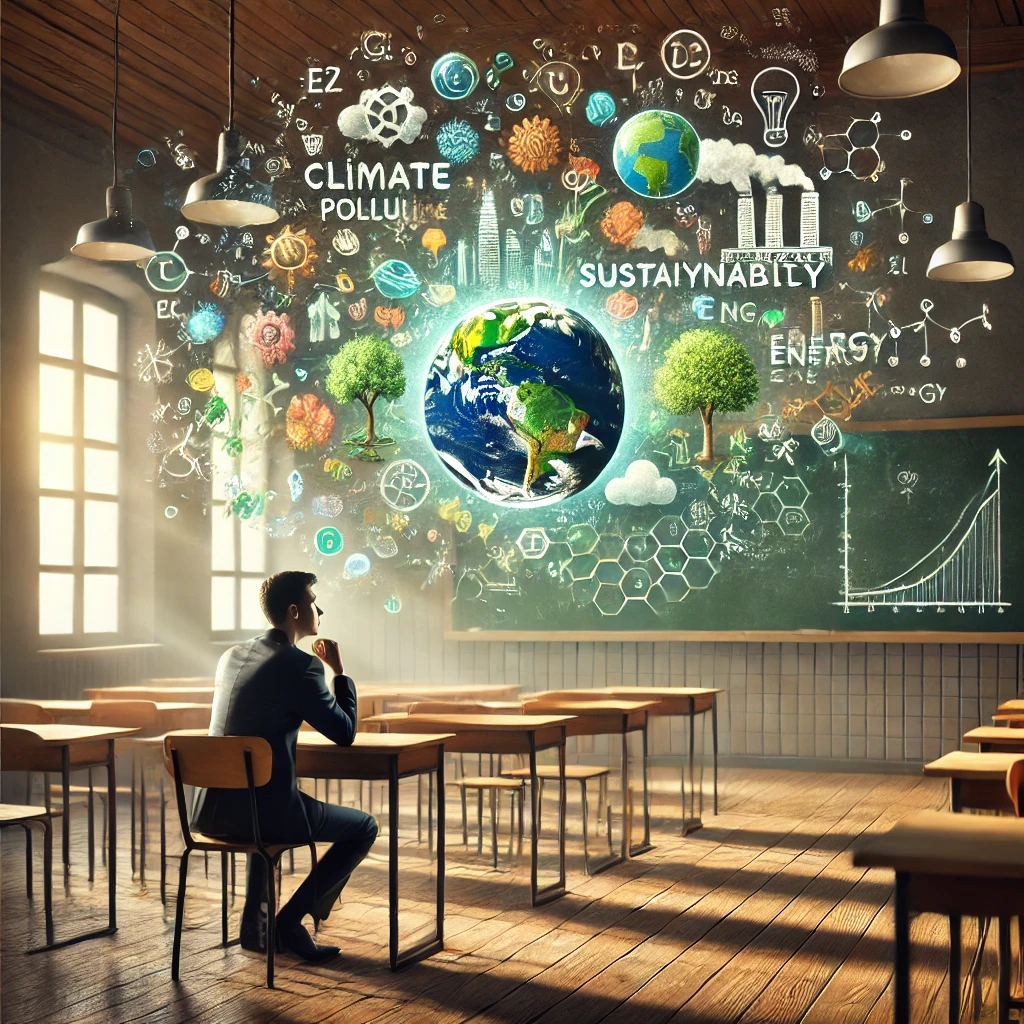}
    \caption*{\bf Sustainability Wanderer    
    }
\end{wrapfigure}

In December 2021, an inquiry was made by the head of master's program in Mathematics at KTH, to 
deliver a course (in Autumn 2022) in applied  mathematics that had components of sustainability. 
After skimming through a  substantial collection of accumulated material, including articles, books, and webpages, totaling over 2 GB of data in 925 files\footnote{The word 'sustainability' appeared over 2.5 billion times in the web search results conducted.} the hardship and challenge of such an undertaking  became more clear.

Starting in early May, the task of reviewing all the material began. After a week, it became clear that working through everything was overwhelming and impractical. The approach was adjusted to focus on five specific topics for which more substantial material was available: Climate, Ecosystems, Pollution, Energy Systems, and Socio-ecology/Socio-economics.
These topics align with the UN Sustainable Development Goals (SDGs), part of the 2030 Agenda, see \cite{sdgs}.

\subsection{Preparing  lecture notes}

The right approach to incorporating mathematics into sustainability  proved to be a more challenging task than initially thought, as there were almost no books or teaching materials—at least none that could be found—specifically combining these two fields at a graduate level. Most of the literature encountered seemed to focus either on the mathematical aspects or on the philosophical dimensions of sustainability, often addressing them in isolation. Although a few sources combined the two in a meaningful way, they were still insufficient for the intended purpose.

This required developing a suitable framework by synthesizing advanced environmental science literature \cite{GM21, HY2013}, including ecosystems, earth sciences, and biology. Starting with these sources and integrating practical applications provided an effective approach to sustainability for master’s-level math students.

One particularly compelling definition encountered during the research was both thought-provoking and informative:\footnote{The simplest way to explain sustainability to a layperson is to appeal to 'common sense' in daily actions.} 
\begin{center}
\fbox{\it Sustainability is best defined by its absence.} 
\end{center}

The lecture notes were carefully reviewed to balance mathematics and sustainability, with revisions made for clarity and precision. Four topics were covered, with some exercises included, though more were still needed, focusing on sustainability, modeling, and mathematics.

The plan included an introductory two-hour lecture to outline the course, meet participants, and adjust if needed, followed by four two-hour lectures each week. Hands-on exercise sessions were also planned, with assistance from a faculty member with expertise in applied mathematics and socio-ecology, who supported the sessions effectively.

\section{Material Covered}

\subsection{The Philosophy}

Since the goal of the course was to raise awareness of sustainability among mathematics students, the first lecture should begin with a general overview of the topic and its philosophy. Therefore, the 17 UN Sustainable Development Goals (SDGs) and some of their targets would be introduced during the lecture.
To set clear expectations, it was emphasized that this course was neither a mathematics course nor a sustainability course. Instead, it aimed to 
{\it  raise awareness of sustainability} and offer a  {\it smörgåsbord of mathematical questions} related to sustainability.

A formal and likely accurate definition of sustainability is as follows:

\fbox{\it 
Sustainability refers to the ability to maintain or support a process continuously over time.}

\begin{wrapfigure}[13]{r}{0.41\textwidth} 
    \centering
    \includegraphics[width=0.32\textwidth]{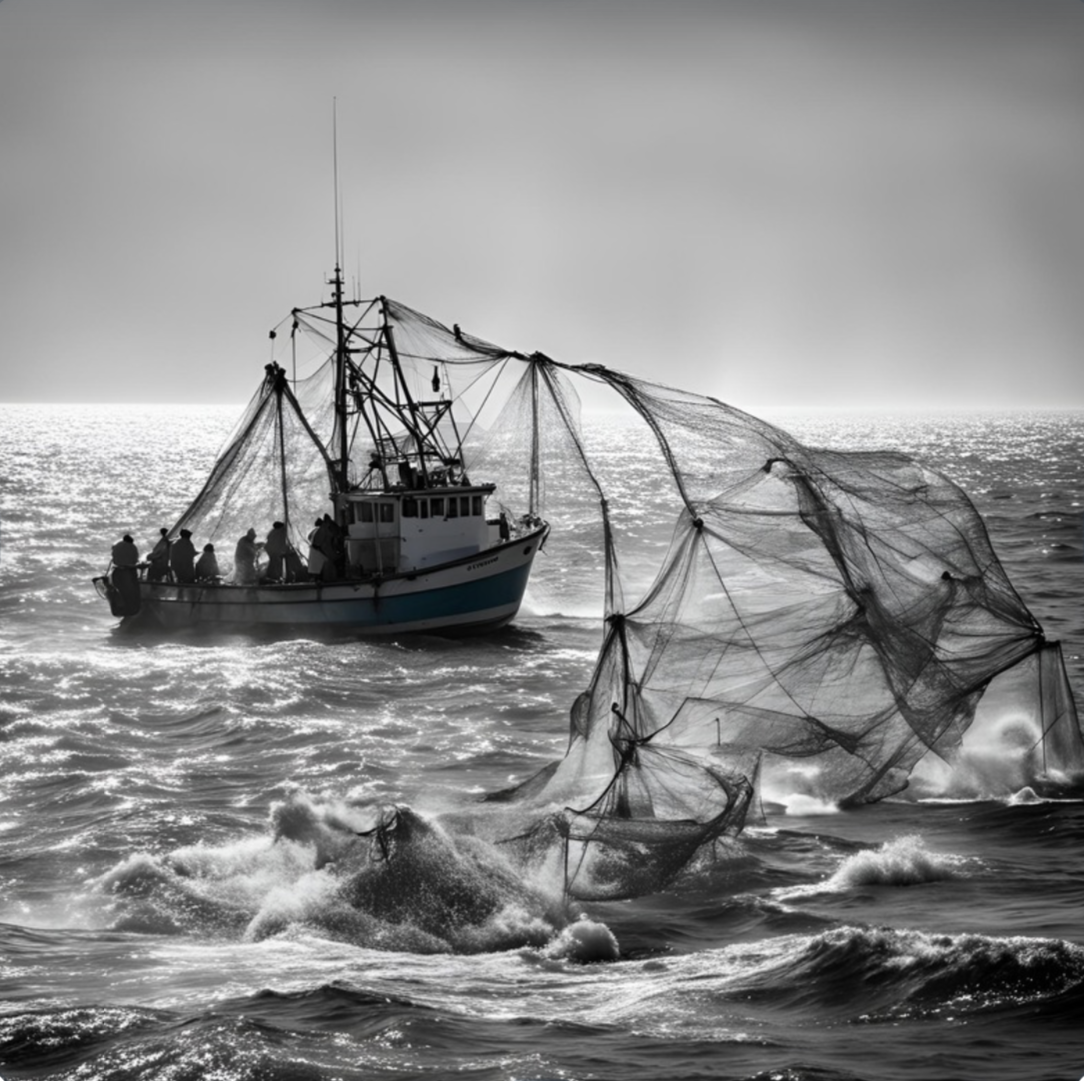}
    \caption*{\bf Overfishing }
    \label{fig:overfishing}
\end{wrapfigure}

This shift in perception shapes our future by aiming to conserve vital resources for future generations. However, sustainability needs contextualization, as it is too abstract. The best way to define it is through examples of unsustainable activities.

The first example  is {\it secondary poisoning}:\footnote{The word {\it relay toxicity} is used in expert discourse and  literature.}
Rat poison not only kills rats but also other  animals eating the poisoned or dead rats. 
For example, California has faced a wildlife crisis involving secondary poisoning from rat poisons for decades, see \cite{urban}.

One may next consider  {\it overharvesting} 
 such as deforestation or overfishing, often occurs at rates too high for species to recover their populations.
The 1990s Atlantic cod collapse revealed the severe ecological and socio-economic impacts of overfishing, displacing over 35,000 workers, see \cite{cod}.

Another example is {\it soil degradation}, which certain agricultural practices, such as slash-and-burn farming, often exacerbate, leading to erosion. In Madagascar  such practices have degraded large amount  of cultivable land, see \cite{madagascar}.

A central theme in the course should be a strong connection to the UN SDGs and relevant targets. Several goals are suitable for mathematical modeling and will be grouped by the mathematical tools used. A brief description of these goals is provided, with further details available at the UN SDGs website.

\begin{wrapfigure}[14]{r}{0.4\textwidth} 
    \centering
    \includegraphics[width=0.35\textwidth]{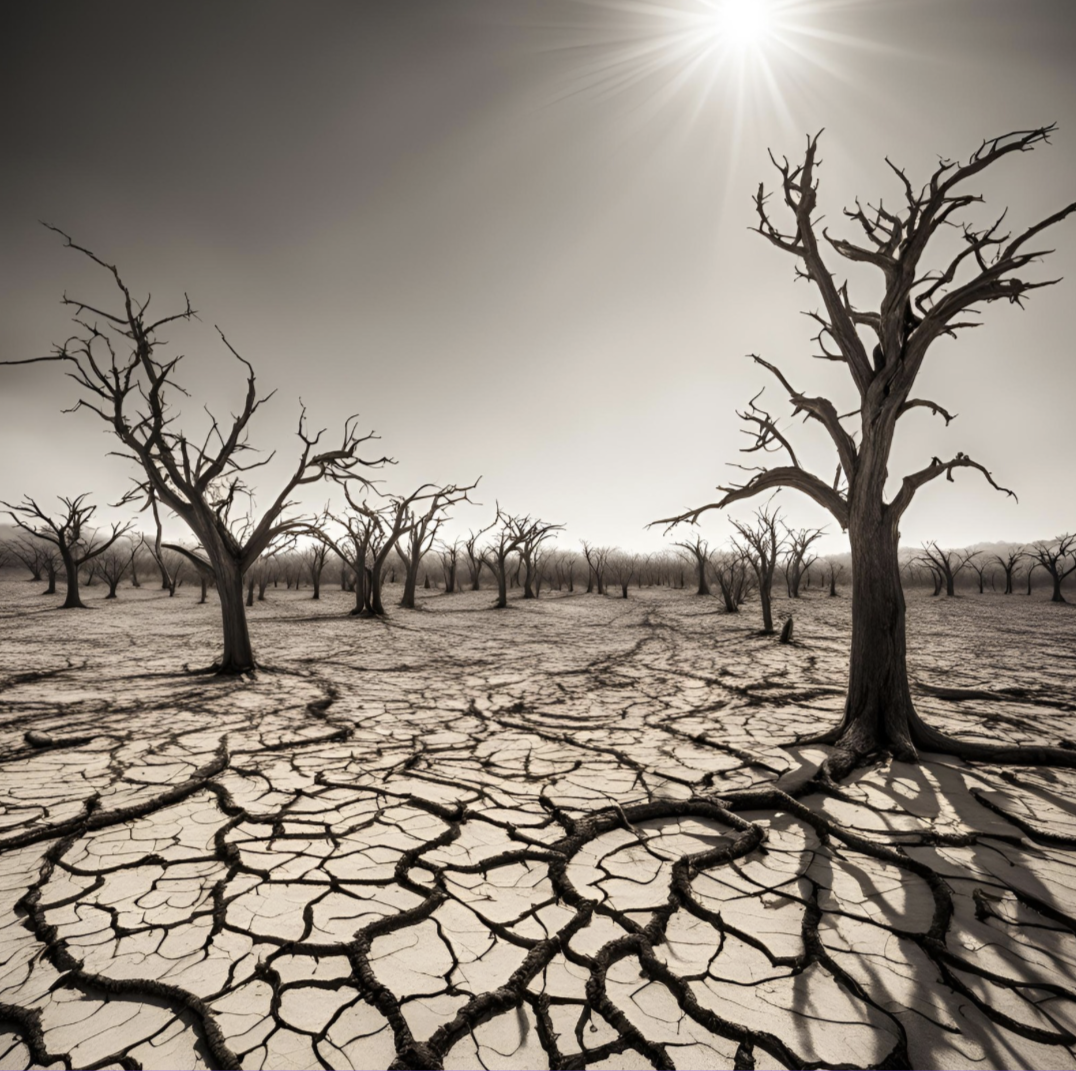}
    \caption*{\bf Soil degradation impacts}
    \label{fig:overfishing}
\end{wrapfigure}

\begin{itemize}

  \item[1)] Goals 7 and 9 relate to sustainable energy and its efficiency, calling for modern energy access for all and the development of resilient infrastructure.

   \item[2)] Goals 4, 5, and 8 relate to health, equality, and socio-economics, calling for inclusiveness and equality, empowering all women and girls, and fostering sustainable economic growth for everyone.
       
   \item[3)] Goals 13, 14, and 15 relate to climate, life on water, and life on land. They call for combating climate change, ensuring the sustainable use of marine resources, and restoring and promoting terrestrial ecosystems.

\end{itemize}

Concerning  goals 7), 9), using mathematical models, such as optimal switching,\footnote{\tiny A simple example is switching between gas, oil, coal, solar, etc., to heat large cities or industrial complexes. } control theory,\footnote{\tiny Control theory uses feedback to optimize wind turbine efficiency in real-time. } and mean field game theory,\footnote{\tiny A mathematical framework for analyzing decision-making in systems with a large number of interacting agents.} 
 can be used to optimize and compare various energy modes, their efficiency, and consumption.
This optimization, subject to constraints, can utilize both partial differential equations and stochastic methods (e.g., to model infrastructure changes and shift analysis).

Goals 4, 5, 8, one can use  several  mathematical models  developed recently, to study  crime, riots, and various other  social/group behaviors directly  linked to social injustices and economic growth.
These models are akin to  {flow dynamics and diffusion.} Similar  {models are associated with international conflicts} and arms races, where behaviors are influenced by  {mass psychology} and the  {impact of governments} through mass media.
These complex models can be examined in their simplest form through  {diffusion equations} (motion under influence).

Concerning goals  13, 14, 15, mathematical  modeling of ecosystems and biodiversity is a  well-established area.
There is  extensive literature  on the study of  carbon dioxide circulation and storage in the ocean. Mathematical models in this field include   thermodynamics, heat flow, the Lorenz equation, precipitation analysis, probability, and ordinary differential equations.

\begin{wrapfigure}[8]{r}{0.4\textwidth} 
    \centering
    \includegraphics[width=0.38\textwidth]{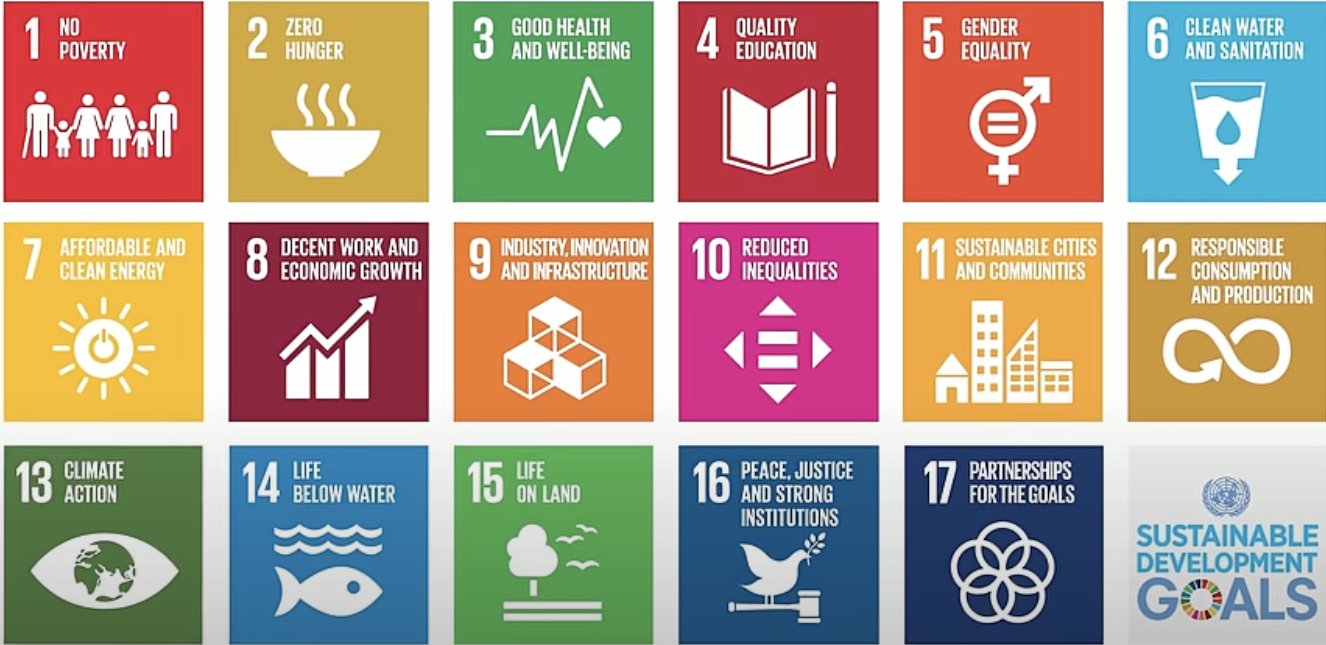}
    \caption*{\bf UN SDG:  {\tiny   \url{https://sdgs.un.org/goals}}
     } 
\end{wrapfigure}

Mathematical modeling is crucial for understanding and leveraging phenomena in nature, industry, biology, economics, and other fields. Technological advancements, driven by mathematics and the sciences, significantly enhance and extend human life. Consequently, there are numerous applications of mathematical thinking and modeling within the field of sustainability.

\vspace{6mm}

To be more precise, one can mention the following areas, among many others:
\begin{itemize}
\item Climate modeling and prediction: PDEs/ODE and numerical simulations.
\item Water resource management: PDEs, and optimizations.
\item Sustainable urban planning: optimization and mean field game theory.
\item Energy systems optimization:  linear programing, and even variational inequalites.
\end{itemize}

Sustainability presents complex challenges, including:
\begin{itemize}
    \item System Complexity: Interactions among environmental, economic, and social systems make comprehensive modeling difficult.
    \item Data Issues: Limited or unreliable data restricts accuracy.
    \item Nonlinear Dynamics/Feedback loops: Small changes can lead to unpredictable, large-scale outcomes.
    \item Collaboration Gaps: Communication barriers and misaligned goals hinder multidisciplinary teamwork.
\end{itemize}

Despite these hurdles, sustainability offers opportunities for innovation, with mathematical modeling playing a key role in addressing these complexities and fostering effective, long-term solutions.

\subsection{Climate}

Understanding Earth's energy budget, the balance between incoming solar radiation and outgoing thermal radiation, is fundamental to studying climate dynamics, addressing climate change, and guiding sustainable practices. 

Climate is characterized by atmospheric variables such as temperature, precipitation, and wind. Its components include radiation, air pressure, humidity, wind patterns, and precipitation. Modeling these parameters is challenging due to the vast range of spatial and temporal scales involved, spanning from seconds to billions of years and from local turbulence to global climate evolution.

Earth's climate system is influenced by interconnected spheres (Atmosphere,\footnote{Governs temperature and weather patterns.}
     Hydrosphere,\footnote{Includes all liquid water on Earth, essential for climate and life.}
     Cryosphere,\footnote{Reflects solar radiation, impacting global temperatures.}
     Lithosphere,\footnote{Forms the physical foundation for continents and oceans.}
     Biosphere\footnote{Encompasses living organisms and their interactions with other spheres.})
and its  modeling  rely on systems of coupled PDEs to describe changes in state variables like temperature, pressure, and velocity.  At the same time, autonomous systems,  allow the use of dynamical systems theory to study the evolution of these models over time.

The basic equations governing the atmosphere involve seven key variables: the velocity components $(u, v, w)$, pressure $p$, temperature $T$, specific humidity $q$, and density $\rho$. These equations describe atmospheric behavior 
and relies on     Newton's Second Law 
    $$
    \frac{d {\bf V} }{dt} = - \frac{1}{\rho}\nabla p - {\bf g} + {\bf F}_{fric} - 2 {\boldsymbol \Omega} \times {\bf V}.
    $$
    Here, $\frac{d}{dt}$ is the material derivative, ${\boldsymbol \Omega}$ is Earth's angular velocity, and ${\boldsymbol \Omega} \times {\bf V}$ is the Coriolis force.
One has to add additional equations to this:
    \begin{itemize}
        \item continuity equation: $\partial_t \rho = - \nabla \cdot (\rho {\bf v})$,
        \item conservation of water vapor: $\partial_t(\rho q) = - \nabla \cdot (\rho q {\bf v }) + \rho (E - C)$,
        \item first law of thermodynamics: $Q = C_p \partial_t T - (1/\rho) dp/dt$,
        \item equation of state: $p = \rho R_g T$.
    \end{itemize}

The response of the climate system to perturbations, such as changes in solar radiation, atmospheric composition, and surface topography, further influences its behavior. While PDEs capture local and detailed processes like heat transfer and fluid flow, dynamical systems provide insights into global, long-term trends and patterns.

The ground work for dynamical system approach was laid by
Edward Lorenz,  who in his investigation of long-term weather forecasting, modeled the Earth's atmosphere as an incompressible fluid situated between two horizontal planes. The fluid is heated from below and cooled at the top, creating a convective flow known as Rayleigh–Bénard flow. This flow is driven by temperature gradients and often forms a regular pattern of cells.

In this simplified model, the flow is reduced to a two-dimensional system with horizontal and vertical velocity components. The boundaries constrain the flow to remain horizontal.

The Lorenz system is a mathematical model that captures the essential dynamics of convective flow. It is described by the following equations:
$$
\dot{x} = \sigma (y - x), \quad \dot{y} = x (\rho - z) - y, \quad \dot{z} = xy - \beta z ,
$$
where $\sigma, \rho, \beta$ are certain parameters and $x$ represents the  convection,
$y$ the  horisontal velocity, and $z$ the temperature difference.

The Lorenz equations are derived by simplifying the complex equations governing atmospheric fluid dynamics. 
 The derivation relies on several key assumptions:
\begin{itemize}
    \item Treating the atmosphere as a fluid.
    \item Reducing the system to three primary variables:
    \begin{itemize}
        \item Convection intensity.
        \item Velocity.
        \item Temperature difference.
    \end{itemize}
    \item Describing the interaction of these variables, which leads to chaotic behavior.
\end{itemize}
It also shows  how small changes in initial conditions can result in unpredictable outcomes, a hallmark of chaotic dynamics. Despite its simplicity, the Lorenz model is a powerful tool for understanding the nonlinear and chaotic nature of convective flows in the atmosphere.

\subsection{Ecosystem}

The study of ecosystems is a branch of ecology, which examines how organisms interact with their physical environment. Ecology is often regarded as a holistic science that treats systems as interconnected wholes rather than isolated parts.

\vspace{-3mm}

\begin{center}
\fbox{An ecosystem consists of a group of organisms together with their  physical environment. }
\end{center}

Ecosystems exhibit intricate interdependencies, making their understanding crucial for addressing global environmental challenges. These systems vary widely in size and diversity, ranging from a small pond to the vast Amazon rainforest, which supports an immense variety of species.

Common examples of ecosystems include agroecosystems, aquatic ecosystems, coral reefs, deserts, forests, human ecosystems, marine ecosystems, rainforests, savannas, and urban ecosystems. The essential components of any ecosystem are water and its temperature, plants, animals, air, light, and soil. Several principles govern ecosystems and their dynamics, and central  concepts include:   {\it geometric growth, cooperation, competition, interacting species, limiting factors},\footnote{Although not well-spelled in the literature, "limiting factors" represent a more complex nonlinearity, and may lead to some interesting variational inequalities once spatial effects are introduced.} 
and many more.

These factors can be modeled mathematically using nonlinear predator-prey models based on principles of population ecology. A simple model of a basic ecosystem consists of a system of ODEs (spatial dynamics can be included using a diffusion term):
$$
\frac{dN_1}{dt} = rN_1 - cN_1N_2, \qquad \frac{dN_2}{dt} = bN_1N_2 - mN_2,
$$
where $r$, $c$, $b$, and $m$ are parameters representing growth rates, interaction rates, and mortality rates. 
However, a  more complex system s given by 
$$
\begin{aligned}
\frac{dx}{dt} & = r_x x \min\left(\frac{R_1}{K_1}, \frac{R_2}{K_2}, \ldots, \frac{R_n}{K_n}\right) - g(x, y, R_i) + \xi_x(t), \\
\frac{dy}{dt} & = \delta g(x, y, R_i) - m y + \xi_y(t), \\
\frac{dR_i}{dt} & = S_i - \phi_i(x, y, R_i), \quad i = 1, 2, \ldots, n,
\end{aligned}
$$
where $x$ represents the prey population and $y$ the predator population. The variables $R_i(t)$ denote the availability of resource $i$ at time $t$. The parameter $r_x$ is the intrinsic growth rate of the prey, while $g(x, y, R_i)$ is the predation rate, which may depend on prey, predators, and resource availability. The conversion efficiency of prey biomass into predator biomass is given by $\delta$, and $m$ is the natural mortality rate of the predator. Resource dynamics are described by $S_i$, the supply rate of resource $i$ (e.g., regeneration or external input), and $\phi_i(x, y, R_i)$, the depletion rate of resource $i$ due to consumption or environmental loss. The term $\min\left(\frac{R_1}{K_1}, \frac{R_2}{K_2}, \ldots, \frac{R_n}{K_n}\right)$ represents Liebig's Law of the Minimum, ensuring that prey growth is constrained by the scarcest resource, where $K_i$ is the maximum utilization level of resource $i$. Finally, the terms $\xi_x(t)$ and $\xi_y(t)$ represent stochastic environmental variability affecting the prey and predator populations, respectively.

A further generalization of the mathematical model arises when considering the so-called meta-ecosystems, which involve interconnected communities, such as river networks where resource flows occur between upstream and downstream ecosystems or between terrestrial and aquatic systems. Key properties of these systems include movement and physical connectivity, which can be represented by a system of PDEs with a connectivity matrix.

\subsection{Pollution}

Pollution is the introduction of harmful substances, known as pollutants, into the environment, degrading the quality of air, water, and land. Pollutants are broadly categorized as either natural (such as volcanic ash, sea salt particles, and photochemically formed ozone) or human-created (like trash and industrial runoff). Furthermore, pollutants can be classified by their level of impact: primary pollutants are emitted directly from sources; secondary pollutants form when primary pollutants react within the environment (for example, smog); and tertiary pollutants represent long-term effects, such as the health consequences of climate change.

All living organisms depend on clean air and water, making pollution a global problem. Urban areas tend to experience higher pollution levels, but it can spread to remote, uninhabited regions through air and water currents. 

\begin{center}
\fbox{\it Pollution threatens all forms of life.}
\end{center}

A major challenge with modelling pollution is that it  can be transported far from its source through air, and water currents.
For instance, accidental releases from nuclear reactors, dispersed worldwide by winds, or smoke from a factory drifting across national borders, highlight the transnational nature of pollution.

Air pollution is a major global health issue and one of the leading risk factors for death, responsible for millions of fatalities annually.\footnote{Data taken from \url{https://ourworldindata.org/air-pollution}} 
It accounts for 11.65\% of global deaths, with the highest burden in low-to-middle-income countries. Premature deaths due to air pollution are projected to increase from around 3 million in 2010 to between 6 and 9 million per year by 2060. Progress in reducing indoor pollution has helped lower global death rates, highlighting the effectiveness of targeted interventions. Beyond health impacts, air pollution has significant economic consequences, including reduced labor productivity, higher health care costs, and lower agricultural yields. Without stringent policies, rising economic activity and energy demand are expected to drive further emissions growth, with economic losses projected to reach up to 1\% of global GDP by 2060.

\medskip
\noindent

To understand the dynamics of noxious fumes transport, it is essential to consider how pollutants are dispersed in the atmosphere. The movement and spread of these fumes are governed by two key processes: advection and diffusion. Advection refers to the transport of pollutants by wind, which carries the fumes while largely maintaining their concentrated profile. Diffusion, on the other hand, involves the spreading of pollutants from areas of high concentration to areas of lower concentration, leading to their gradual dispersion. Together, these processes determine the transport, distribution, and eventual dilution of noxious fumes in the environment.

To account for both advection and diffusion, the governing equation is given by
\[
\frac{\partial C}{\partial t} + \vec{v} \cdot \nabla C = D \nabla^2 C + S - R + Q(C),
\]
which describes how concentration \( C \) evolves, including chemical reactions \( Q(C) \), spatial diffusion, and source \( S \) and removal \( R \) terms. 

Emission control and regulation models are essential for managing air pollution within complex physical, economic, and regulatory systems. A linear emission reduction model, often a first approximation, expresses emissions after control as
\[
E_c = E_0 (1 - R),
\]
where \( E_0 \) is the baseline emissions, \( R \) is the reduction rate, and \( E_c \leq E_{\text{limit}} \) ensures compliance with regulatory caps. In practice, emission reduction costs are nonlinear, modeled using marginal abatement cost functions \( C_i(x_i) \), with an optimization problem:
\[
\text{Minimize} \, C = \sum_{i=1}^n C_i(x_i), \quad \text{subject to} \quad \sum_{i=1}^n (E_i - x_i) \leq E_{\text{limit}},
\]
where \( x_i \) is the reduction at source \( i \).

\subsection{Energy}

Energy is the cornerstone of modern society, providing the foundation for global infrastructure. However, the world is now facing its first  global energy crisis, driven by overconsumption, poor infrastructure, geopolitical conflicts, and natural disruptions. The solution lies not in further reliance on fossil fuels but in advancing the energy transition through coordinated global efforts and a focus on clean, sustainable resources.

The future of energy is shaped by factors such as policy-making, rising decarbonization costs, energy supply challenges, and government-driven "greenflation" (inflation caused by increasing costs of green technologies). Ensuring energy security requires integrating various renewable sources:  solar energy, the most abundant resource; wind energy, enhanced by taller turbines and larger rotors; geothermal energy, which harnesses Earth’s internal heat; hydropower, relying on water movement; ocean energy, derived from seawater’s properties; and bioenergy, sourced from organic materials like wood and crops.

Several key factors drive the energy transition, including decreasing renewable technology costs, the need for modern grid infrastructure, tax incentives, subsidies, and a skilled workforce. Despite these advancements, the transition is financially demanding, often exceeding the capacities of companies and nations. International collaboration is essential to achieve the United Nations' Sustainable Development Goals  related to energy.

Mathematical and computational models are critical for navigating this transition, combining energy systems engineering, economics, environmental science, and behavioral studies. Tools like optimization theory, game theory, and numerical analysis help assess energy scenarios, economic impacts, and sustainability outcomes. 
Important  metrics that guide these evaluations are:

\begin{itemize}
    \item Energy Return on Energy Investment (EROEI): The ratio of energy produced to energy invested.  
    \item Energy Payback Time (EPBT): The time required for a system to generate the energy used in its creation.  
    \item Levelized Cost of Energy (LCOE): The average cost of energy production over a system’s lifetime.  
    \item Capacity Factor (CF): The ratio of actual output to maximum potential output.  
    \item Carbon Intensity: Greenhouse gas emissions per unit of energy produced.  
\end{itemize}

EROEI is defined as the ratio of the amount of final usable energy acquired from a particular energy resource to the amount of primary energy expended to obtain that energy resource:
$$
\text{EROEI} = \frac{\text{Usable Acquired Energy}}{\text{Energy Expended to Get that Energy}}.
$$
More specifically, this ratio—typically dimensionless—means that a given energy production technology will provide EROEI joules (J) on an energy investment of 1 J.
Alternatively one can define  the continuous EROEI model, which uses ordinary differential equations (ODEs):

$$
    \frac{dR}{dt} = -\frac{E_o(t)}{\eta}, \qquad  
    E_o(t) = \kappa R(t)^n,  \qquad  
    E_i(t) = E_{i0} + \beta (R_{\text{max}} - R(t)),  \\ 
    $$
   $$
    \text{EROEI}(t) = \frac{\kappa R(t)^n}{E_{i0} + \beta (R_{\text{max}} - R(t))}.
    $$

Here, \( R(t) \) represents the remaining energy resource at time \( t \), \( E_o(t) \) is the energy output at time \( t \), \( E_i(t) \) denotes the energy investment at time \( t \), and \( \text{EROEI}(t) = \frac{E_o(t)}{E_i(t)} \) is the Energy Return on Investment.\footnote{The parameters used in the model are as follows: \( \eta \) represents the extraction efficiency, \( \kappa \) is the proportionality constant for output, \( n \) denotes the nonlinearity of output with remaining resources, \( \beta \) is the investment scaling factor, \( E_{i0} \) refers to the initial energy investment, and \( R_{\text{max}} \) is the initial total resource.
}

\subsection{Tipping Point}

\begin{wrapfigure}[12]{r}{0.4\textwidth} 
    \centering
    \includegraphics[width=0.32\textwidth]{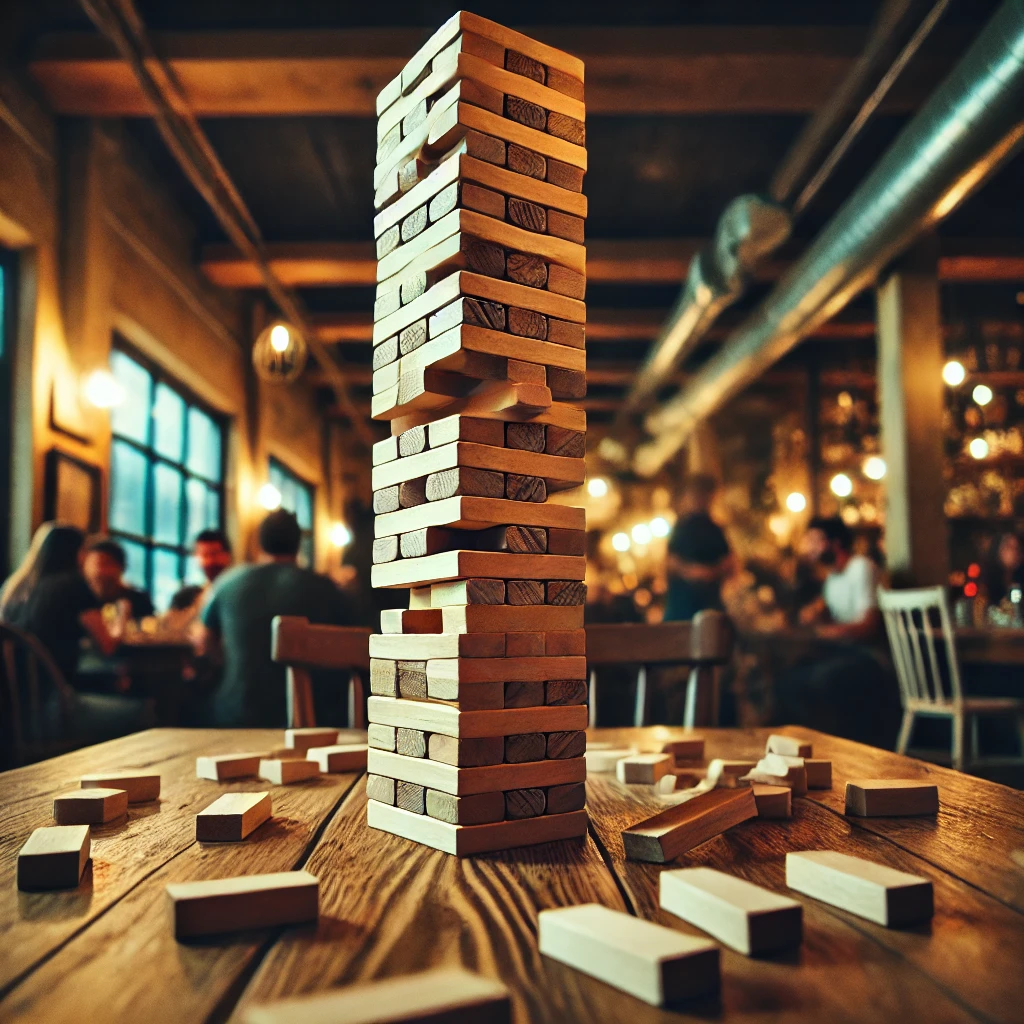}
    \caption*{\bf  Jenga Game}
\end{wrapfigure}

A tipping point refers to a critical threshold at which a small change in input or conditions can lead to a significant, often irreversible change in the state of a system.
Tipping points are seen across various disciplines, including climate science, biology, social science, and economics.

In the game Jenga, the goal is to remove and stack blocks while keeping the tower stable. The tipping point occurs when the system becomes too unstable, causing the tower to collapse—a change that cannot self-correct.

\bigskip

We give a few examples of tipping points:

\begin{itemize}
\item The mountain mist frog:  Native to northern Queensland's wet forests, this species  was declared extinct by the IUCN in 2022 after disappearing in the 1990s.
    \item Arctic Sea Ice Loss: A critical threshold where warming causes rapid and potentially irreversible melting of sea ice, reducing Earth's albedo and accelerating climate change.
    \item Amazon Rainforest Dieback: A tipping point where deforestation and climate change trigger a feedback loop, converting rainforest into savannah and releasing massive carbon stores.
    \item Coral Reef Bleaching: Rising sea temperatures lead to mass coral bleaching events, potentially pushing reefs beyond recovery, drastically impacting marine ecosystems.
    \item Greenland Ice Sheet Melting: Crossing a threshold where ice loss becomes self-reinforcing, leading to significant sea-level rise and further global warming.
\end{itemize}

Sweden’s rising gang violence (2022–2023) marks a critical tipping point in urban dynamics. Socio-economic divides, weak integration policies, and escalating organized crime have fueled a surge in shootings and bombings, severely straining public trust and governance. Addressing this systemic crisis demands holistic reforms to prevent further escalation and restore stability.

A mathematical model describing crime rate dynamics ($ x(t) $) in an urban setting, 
can be modeled by the equation
$$
\dot{x}(t) = a x(t) - b T(t) - c u_2(t) + d u_1(t).
$$
The model includes socioeconomic factors $u_1(t)$ (unemployment levels), police presence $u_2(t)$ (a controllable factor to reduce crime), and community trust $T(t)$ (representing social cohesion). 
The parameters $a, b, c$, and $d$ respectively represent the natural crime growth rate, the crime reduction sensitivity to community trust, the crime reduction sensitivity to police presence, and the crime increase sensitivity to unemployment.

A more complex model of  dynamics of  control strategy for police deployment, minimizing crime through adaptive responses is given by:
$$
\frac{\partial C}{\partial t} = D \nabla^2 C + R(C) - P(x, t),
$$
where $ C(x, t) $ is crime density at spatial location $ x $ and time $ t $, and
 $D \nabla^2 C$, representing the spatial spread of crime, a reaction term $R(C)$ for intrinsic crime growth, and $P(x, t)$, which accounts for police intervention to reduce crime density.
The police force is distributed dynamically in space and time to solve an optimal control problem:

$$
\min_{P(x, t)} \int_0^T \int_\Omega \left[ C(x, t)^2 + \alpha P(x, t)^2 \right] dx \, dt.
$$
Here, $\alpha$ is the weighting parameter for police deployment costs, $\Omega$ represents the spatial domain, and $T$ denotes the time horizon.
The constraint on police resources is given by 
$$\int_\Omega P(x, t) , dx \leq P_{\text{max}},$$
 illustrating how tipping points in crime emerge and emphasizing adaptive resource allocation to prevent crises and maintain stability.

\newpage

 \section{Curriculum Delivery}
 



\subsection{Administrative support}

The Sustainability Office at our university 
spearheads the institution's commitment to integrating environment and sustainable development across education, research, and collaboration, as dictated by its policy to be a leader in this field. Responsible for the ISO 14001 certified environmental management system,\footnote{https://www.iso.org/standard/60857.html}
 this office works centrally with leadership to minimize the university's environmental impact through operational practices and responsible resource consumption, while also promoting social responsibility. Critically, through this robust policy framework, strong support is provided for the development of courses across various disciplines that integrate articulated components of sustainability, ensuring the institution actively contributes to sustainable development goals through its core academic and collaborative activities. Furthermore, this policy has directly facilitated the development of the current course, offering significant financial and logistical support, which has been instrumental in its development, thereby demonstrating the institution’s commitment to translating its sustainability policy into tangible educational opportunities.

\subsection{Students Background \& Prerequisites }

The course  is  designed for first-year master's students in both pure and applied Mathematics, who meet all formal prerequisites, including a strong foundation in basic mathematics at the bachelor's level, and  the completion of a bachelor's degree project.\footnote{We also have a technology major requirement, which can be relaxed.} 
However, the course can also be adapted for students with less mathematical background, such as Ph.D. students from other disciplines, by modifying the homework assignments and modeling problems. This, however, would require significantly more effort.

\subsection{Goals and Outcomes}

The course aims to equip students with a strong foundation in research methodology, emphasizing its relevance to their respective subjects and professional roles. Through a series of modular learning activities, primarily conducted as seminars, students will engage with topics related to gender equality, diversity, and sustainable development, with a particular focus on achieving a climate-neutral society.

Upon successful completion of the course, students will:

\begin{itemize}
    \item Develop a comprehensive understanding of the role of science and technology in society and the ethical responsibility of individuals in applying their knowledge for sustainable development.  
    \item Effectively communicate and justify their conclusions, demonstrating clarity in presenting their knowledge and arguments both orally and in writing, while engaging in meaningful dialogue with diverse audiences.  
    \item Gain insight into contemporary research and advancements in the field of mathematics, fostering a deeper appreciation for ongoing developments in the discipline.  
\end{itemize}

\subsection{Lectures}

The course consisted of 2–3 hours of lectures per week, including exercise sessions. A teaching assistant with expertise in both mathematics and its applications played a key role. This allowed the exercise sessions to focus primarily on mathematical concepts and problem-solving.

Beamer presentations proved to be extremely helpful. This format allowed for the inclusion of footnotes in smaller fonts, which provided explanations for words and terminology that were unclear or had been forgotten.
Throughout the classes, valuable insights were gained from the students' questions, highlighting areas for improvement in future iterations of the course.

The course included bi-weekly homework assignments that covered both mathematical and sustainability aspects. The latter often required explanations and examples of specific concepts. For instance, students were asked to concisely explain Earth's energy budget on a single page or provide examples of unsustainable activities, meta-ecosystems and so on.

\subsection{Examination}

The examination process for this course consists of multiple components, including homework, a final report, and a presentation. Students were asked to form groups of 4–5 members and select a topic at a very early stage, possibly by the second week. The groups were monitored throughout their work and received guidance from the instructor.

In terms of ethics, all members of a group are responsible for the group's work. During any assessment, every student is required to honestly disclose any help received and sources used. Furthermore, in an oral assessment, every student must be able to present and answer questions about the entire assignment and solution, demonstrating a comprehensive understanding of the material.

For the final exam, students were required to create a Beamer presentation on a course-related topic of their choice. A two-week break from classes was given for preparation, and students had the opportunity to consult with the instructor during office hours. This format provided a positive and effective way for students to learn independently and develop their presentation skills.

Suggested topics for presentation included: Tipping point analysis,  Metrics in sustainability,  Energy modes, Cost of sustainability/transition, 
 ESG: in general framework, and applications, 	
 Socioecology, 
 Socioeconomics,
 Climate and environment,
Traffic flow,
 Complex industrial modelling,
 Biological systems,
 Green Finance,
 Catastrophe theory.

\subsection{Students' Perception and Assessment}  

The first time this course was offered, it was optional for master's students in pure mathematics and had low attendance for various reasons. One key factor appeared to be its less mathematical nature, as most mathematics master's students tend to prioritize high grades and focus on mathematics-intensive courses to strengthen their prospects for third-cycle studies. Another significant factor affecting attendance was a parallel mandatory course on sustainability, which focused on gender equality. 

Following a challenging experience and collaborative discussions with the parallel course instructor, we merged the two courses. The revised curriculum dedicated two weeks to gender equality and intersectional design, with a reduced emphasis on mathematical content, and five weeks to my accumulated modeling material.

Students with a background in pure mathematics are generally less receptive to modeling, particularly when the course includes components beyond traditional mathematics. In this case, the philosophical aspects of sustainability, along with discussions on ethical, moral, and environmental issues, posed an initial challenge to engagement.  

Similarly, students from applied mathematics shared some of this hesitation but showed a greater interest in the modeling aspects of the course. Despite an initial sense of skepticism toward the non-mathematical components, students gradually became more engaged as they encountered a wealth of intriguing information about the environment and related topics. These aspects sparked their curiosity and encouraged active participation in discussions.  

Once the initial phase of the course passed and students began working on hands-on problems in mathematical modeling, their engagement increased significantly. They showed strong interest in exploring complex and thought-provoking questions, leading to high levels of participation in both lectures and discussion-based seminar sessions. This shift in engagement highlighted the importance of gradually integrating mathematical problem-solving with broader interdisciplinary discussions to foster deeper involvement and learning.

\subsection{Challenges}

Throughout the course, several challenges were encountered, particularly in managing group dynamics, ensuring academic integrity, and adapting to evolving learning tools.

One significant challenge was group inhomogeneity, where students had varying levels of mathematical background and problem-solving skills. This affected the balance of workload and collaboration within groups. Additionally, helping students with problems required addressing different levels of understanding, making it essential for the instructor to provide tailored guidance. The instructor’s knowledge of topic variations also played a crucial role, as students often approached problems from different perspectives, requiring flexibility in explanations and problem-solving strategies.

Another ongoing challenge was modifying problems from one semester to the next to prevent repetition and discourage students from relying on previous solutions. Closely related to this was the issue of preventing plagiarism, particularly in an era where students have access to vast online resources.

The increasing use of AI-assisted tools presented both opportunities and challenges. While AI helped students structure their reports and improve coherence, it also led to situations where students worked only on their assigned sections without truly engaging with the work of their peers. As a result, many did not gain a full understanding of the entire project.

To address these issues, we recognized the importance of establishing a strong collaborative spirit within groups at an early stage. Ensuring that all members actively participated was essential to fostering meaningful learning experiences. Additionally, strict monitoring and continuous follow-up on project development became necessary to track progress, encourage engagement, and verify that all students were contributing fairly.

Moving forward, these insights will guide adjustments to the course structure, with a focus on improving group collaboration, academic integrity, and deeper learning outcomes.

\section{Conclusion and Pathways Forward}

Preparing this course allowed for exploring sustainability from a mathematical perspective, further discussed in the conclusion. A key lesson is that broad topics like sustainability are best approached through real-world applications, rather than fitting research into a narrow focus. This approach, more suited for specialized courses, helps identify the right mathematical tools for addressing real-world problems.

The research on the intersection of mathematics and sustainability is still in its infancy, far from being in its formative stages. A comparison can be made with modern mathematical finance, which began to take shape in the mid-20th century when mathematicians and economists started formalizing models for financial markets. One might cite pioneers such as Paul Samuelson (1960s), and Fischer Black, Myron Scholes, and Robert Merton (1973). In contrast, mathematical sustainability is still in its 'prehistory' compared to the development of modern mathematical finance. It will likely take many years before tools that are sufficiently robust to address sustainability-related problems at a more general level are developed.

There are several mathematicians working on applied problems in various fields that directly connect to specific aspects of sustainability, often problem-specific in nature. As a result, it will be very challenging to identify a 'Black-Scholes' equation (see \cite{BS}) 
 for sustainability, or a Mean-Field Game theory (as proposed by Lasry and Lions (see \cite{LL}).
It is also unclear whether there exists a unifying mathematical model that can encompass many, if not most or all, aspects of sustainability-related problems.

In conclusion, the involvement of theoretical mathematicians is crucial for the field to take root, and broader participation will be essential for its continued advancement.
As noted earlier, progress is currently fragmented, with research focusing on specific problems. Whether a more unified theory will emerge remains to be seen.

Several key areas in sustainability offer intriguing opportunities for mathematical exploration, as discussed throughout this text. These include tipping point analysis, sustainability metrics, modeling feedback mechanisms, delay equations, variational formulations incorporating spatial parameters, and the application of Mean Field Game theory. Many other promising directions likely remain undiscovered.

\section{Acknowledgements}
I would like to express my gratitude to several individuals and organizations for their support, encouragement, and discussions. These include Lilit Afram Shahgholian (Swedavia), Katharina Gustafsson (KTH), Fredrik Viklund (KTH), Gunnar Tibert (KTH), Maja Schlüter (Stockholm Resilience Center), Peter Tankov (ENSAE), Sonja Radosavljevic, Ulf Näslund (Vasakronan), Uno Wennergren (Linköping University), and KTH Sustainability Office (for a generous financial support).

\end{document}